\magnification=\magstep1
\input amstex
\documentstyle {amsppt}
\input epsf
\NoBlackBoxes
\NoRunningHeads
\pageheight {8.5 truein}
\pagewidth {6.5 truein}

%PAPER FOR PROCEEDINGS OF THE NC STATE CONFERENCE

\topmatter
\title
Construction of some algebras associated to directed graphs
and related to factorizations of noncommutative polynomials
\endtitle
\author
Vladimir Retakh, Shirlei Serconek  and Robert Lee Wilson 
\endauthor
\address
\newline
V.R., R.W.: Department of Mathematics, Rutgers University,
Piscataway,
NJ 08854-8019
\newline
S.S:
IME-UFG
CX Postal 131
Goiania - GO
CEP 74001-970 Brazil
\endaddress
\email
\newline vretakh\@math.rutgers.edu
\newline
serconek\@math.rutgers.edu
\newline 
rwilson\@math.rutgers.edu
\endemail
\abstract This is a survey of recently published results. 
We introduce and study a wide class of algebras associated to 
directed graphs and related to factorizations of noncommutative 
polynomials. In particular, we show that for many well-known 
graphs such algebras are Koszul and compute their Hilbert series.
\endabstract

\endtopmatter
\document

Let $R$ be an associative ring with unit and
$P(t)=a_0t^n+a_1t^{n-1}+\dots +a_n$ be a polynomial over $R$. Here
$t$ is an independent central variable. We consider
factorizations of $P(t)$ into a product
$$P(t)=a_0(t-y_n)(t-y_{n-1})\dots (t-y_1) \tag 0.1$$
if such factorizations exist.

When $R$ is a (commutative) field, there is at most one such factorization up
to a permutation of factors. When $R$ is not commutative, the
polynomial $P(t)$ may have several essentially different
factorizations.

The set of factorizations of a polynomial over a noncommutative
ring can be rather complicated and studying them is a challenging
and useful problem (see, for example, \cite {N, GLR, GR1, GR2,
GRW, GGRW, LL, O, B, V, W}).

In this paper we present an approach relating such factorizations
to algebras associated with directed graphs and study properties
of such algebras.

In the factorization (0.1) the element $y_1$ is called a {\it right root}
of $P(t)$ and element $y_n$ is called a {\it left root} of
$P(t)$.  This terminology can be justified by the following
equalities (see, for example, \cite {L}):
$$a_0y_1^n+a_1y_1^{n-1}+\dots +a_{n-1}y_1+a_n=0,  $$
$$y_n^na_0+y_n^{n-1}a_1+\dots +y_na_{n-1}+a_n=0.$$

It is natural (see \cite {GGRSW}) to call the elements $y_i$,
$i=1,2,\dots, y_n$, in (0.1) {\it pseudo-roots} of $P(t)$. Note
that any root is a pseudo-root but not every pseudo-root is a root.

There is a connection between factorizations of noncommutative
polynomials and finding their right (left) roots and
pseudo-roots. Consider the simplest non-trivial example when
$P(t)=t^2+a_1t+a_2$ and $x_1, x_2$ are its right roots such that
the difference $x_1-x_2$ is invertible in $R$.

Set $$x_{1,2}=(x_2-x_1)x_2(x_2-x_1)^{-1}, $$
$$x_{2,1}=(x_1-x_2)x_1(x_1-x_2)^{-1}.$$

One can show that
$$P(t)=(t-x_{1,2})(t-x_1)= (t-x_{2,1})(t-x_2). \tag 0.2$$

Thus we have two different factorizations of $P(t)$.

In studying factorizations of noncommutative polynomials, one
should answer at least two questions:

1) How to obtain factorizations of type (0.1)?

2) How to relate different factorizations?

A partial answer to the first question was given in \cite {GR1,
GR2} when a polynomial $P(t)=t^n+a_1t^{n-1}+\dots +a_n$ has $n$
right roots
 $x_1, x_2,\dots , x_n$  in a {\it generic position}.
Here elements $x_1, x_2,\dots , x_n$ are in a generic position if
all Vandermonde matrices $(x_{i_{\ell}}^{m-1})$, $\ell
,m=1,2,\dots , k$ are invertible over $R$ when $k\geq 2$ and
$i_1, i_2,\dots , i_k$ are distinct. In \cite {GR1, GR2} explicit
formulas for $n!$ factorizations of $P(t)$ were found. For $n=2$
these formulas give the factorizations (0.2).

To answer the second question, one has to study relations among
pseudo-roots of $P(t)$ and, more generally, study properties of
subalgebras generated by pseudo-roots of $P(t)$.

For example, the factorizations (0.2) imply the following identities
between the corresponding pseudo-roots:
$$x_{1,2}+x_1= x_{2,1}+x_2 ,\tag 0.3a$$
$$x_{1,2}x_1= x_{2,1}x_2. \tag 0.3b$$

Note that one can consider expressions in (0.3) as noncommutative
elementary symmetric functions of order $1$ and $2$ in $x_1, x_2$
(see \cite {GR1, GR2, GRW, GGRW}). It is natural to study the algebra
$Q_2$ with generators $x_{1,2}, x_{2,1}, x_1, x_2$ satisfying
relations (0.3).

To carry out this idea, a {\it universal algebra of } $n$ {\it
pseudo-roots} in a generic position, called $Q_n$, was introduced
in \cite {GRW} and its properties were studied in detail in
\cite {GRW, GGR, GGRSW, SW, Pi}. The algebras $Q_n$ are defined by
linear and quadratic relations similar to relations (0.3). These
algebras are Koszul  and have nice  Hilbert series. Overall, one
can say that $Q_n$ is a rather ``tame" algebra despite its
exponential growth.

A more general approach to a study of algebras of pseudo-roots
starts with $\Gamma _P$,  {\it a directed graph of right divisors
of a polynomial} $P(t)$ similar to the graph of divisors of a
natural number. In this graph {\it vertices} correspond to right
divisors of $P(t)$ and {\it edges} to pseudo-roots of $P(t)$.
Factorizations (0.1) correspond to {\it paths} from a maximal
vertex $P(t)$ to a minimal vertex $1$.

\midinsert
\centerline{
\epsfxsize=2in
\epsfbox{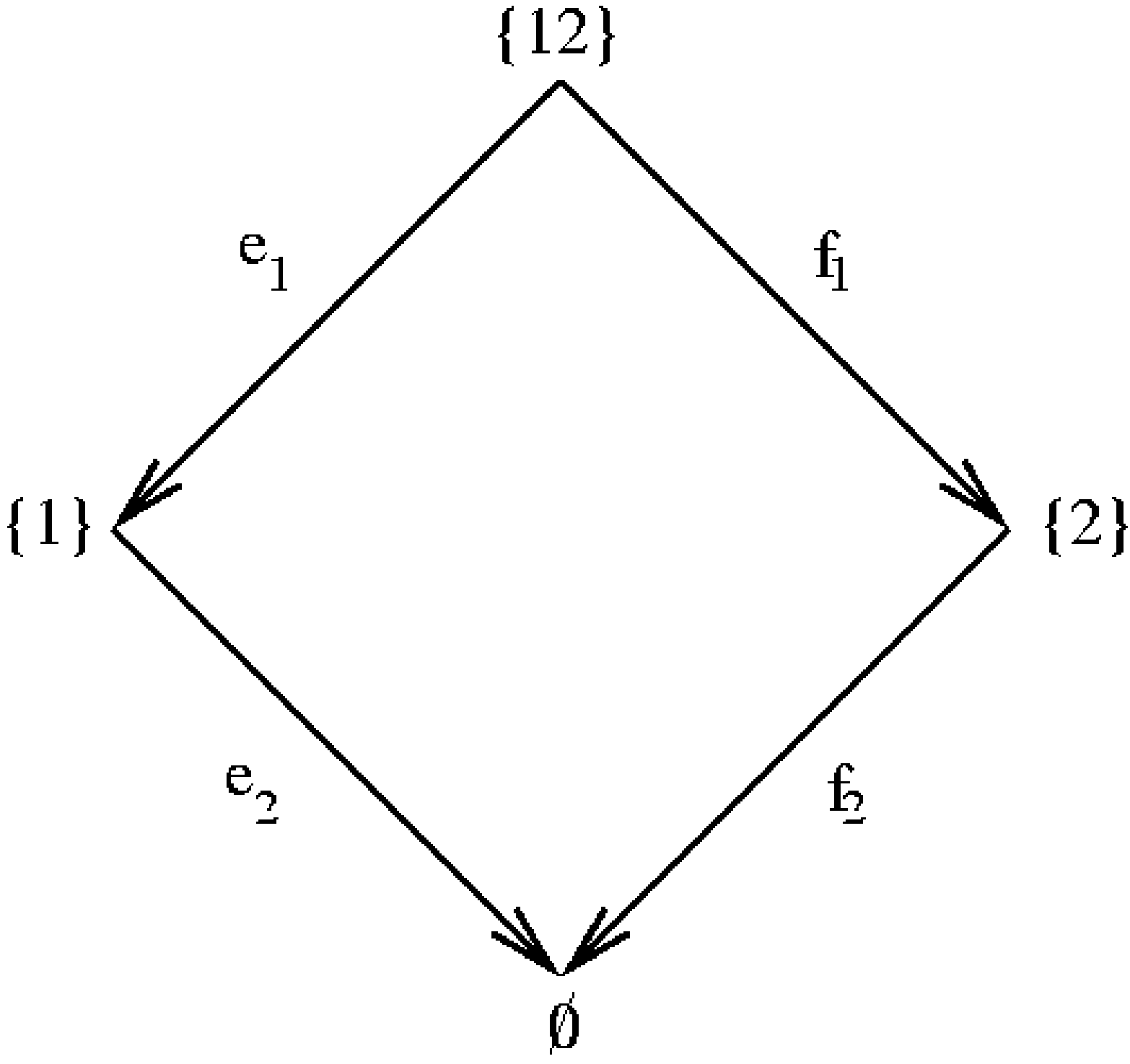}}
\botcaption{Figure 1}
\endcaption
\endinsert

In this setting relations (0.3) for the algebra $Q_2$ can be
described with a help of the {\it diamond} graph with vertices
indexed by subsets of set $\{1,2\}$ (see Fig 1). Elements
$x_{1,2}$ and $x_{2,1}$ correspond to edges $e_1$ and $f_1$,
and elements $x_1$ and $x_2$ correspond to edges $e_2$ and $f_2$
respectively. Factorizations (0.2) correspond to two paths with
the same origin and the end.

Our main objects are the ``universal algebras" $A(\Gamma )$
associated to a directed graph $\Gamma $  and ``universal
polynomials" over these algebras. These algebras are generated by
edges of $\Gamma $ and relations are defined by pairs of directed
paths in $\Gamma $ with the same origin and end.  If $\Gamma
$ is the Hasse graph of the Boolean lattice of subsets of the set
$\{1,2,\dots , n\}$ then $A(\Gamma )=Q_n$.

It turns out that the algebras $A(\Gamma )$ have a lot of interesting
properties similar to the algebras $Q_n$.   The ``geometric nature" of
algebras $A(\Gamma )$ helps us to understand their structure and
to simplify the proofs of main results for algebras $Q_n$ compared
to the proofs given in \cite {GRW, GGRSW, SW, Pi}. A geometric
approach expresses coefficients of the ``universal polynomials"
via certain pseudo-roots of these polynomials, giving a geometric
version of the noncommutative Vi\'ete theorem from \cite {GR1,
GR2}, see also \cite {GGRW}.

This paper contains definitions and results from \cite {GRSW, RSW,
RSW1, GGRW1}. It is organized in the following way.  In Section 1
we recall  the noncommutative Vi\'ete theorem and the construction of
the algebras $Q_n$. Section 2 contains a definition of the algebras
$A(\Gamma )$ and similar algebras. In Section 3 we describe a
linear basis in $A(\Gamma )$ and state that for a large class of
directed graphs, the algebras $A(\Gamma )$ are Koszul. In Section 4 we
compute Hilbert series for algebras $A(\Gamma )$ and discuss some
examples. In Section 5 we describe a geometric version of the
noncommutative Vi\'ete theorem.

During preparation of this paper
Vladimir Retakh  and Robert Lee Wilson were partially supported by NSA, 
Shirlei Serconek was partially supported by CNPq/PADCT.

\head 1. The Vi\'ete Theorem and algebras of pseudo-roots of
noncommutative polynomials\endhead

\subhead 1.1. Factorizations of noncommutative
polynomials\endsubhead
 A generalization of the factorization (0.2)
can be presented as follows (see \cite {GR1, GR2, GRW}). Let $R$
be an associative ring with unit and $P(t)=t^n+a_1t^{n-1}+\dots
+a_n$ be polynomial in $R[t]$. Let $x_1, x_2,\dots , x_n\in R$ be
right roots of the polynomial, i.e. $x_i^n+a_1x_i^{n-1}+\dots
+a_n=0$ for all $i$.

We say that the roots $x_1, x_2,\dots , x_n$ are in {\it generic
position} if all Vandermonde matrices $V(x_{i_1}, x_{i_2},\dots ,
x_{i_k})= (x_{i_l}^{m-1})$, $l,m=1,2,\dots , k$ are invertible over
$R$ when $k\geq 2$ and all $i_1, i_2,\dots , i_k$ are distinct
and every entry of such an inverse matrix is invertible.

In this case, one can define the Vandermonde quasideterminant
(for the general theory of quasideterminants see \cite {GR, GR2,
GGRW})
$$v(x_{i_1}, x_{i_2},\dots , x_{i_k})=
x_{i_k}^{k-1} - r\cdot V(x_{i_1}, x_{i_2},\dots
,x_{i_{k-1}})^{-1}\cdot c
$$
where $r$ is a row vector of length $k-1$, $r=(x_{i_1}^{k-1}, x_{i_2}^{k-1},\dots
, x_{i_{k-1}}^{k-1})$ and $c$ is a column vector of the same length, $c=(1, x_{i_k},
\dots , x_{i_k}^{k-2})^T$.

One can see that $v(x_{i_1}, x_{i_2},\dots , x_{i_k})$ is an
element in the ring $R$. This element is invertible in $R$
because elements $x_1, x_2,\dots, x_n$ are in a generic position.
\smallskip
\noindent {\bf Example}. $v(x_1,x_2)=x_2-x_1$.
\smallskip

Note that $v(x_{i_1}, x_{i_2},\dots , x_{i_k})$ is a rational
expression in $x_{i_1}, x_{i_2},\dots , x_{i_k}$ and that it does
not depend of the ordering of $x_{i_1}, x_{i_2},\dots , x_{i_{k-1}}$
(see \cite {GR, GR2, GGRW}).

Set
$$x_{A,i_k}=v(x_{i_1}, x_{i_2},\dots , x_{i_k})\cdot x_{i_k}\cdot
v(x_{i_1}, x_{i_2},\dots , x_{i_k})^{-1},$$ where $A=\{x_{i_1},
x_{i_2},\dots , x_{i_{k-1}}\}$. Our notation is justified because
$x_{A,i}$ also does not depend on the ordering of $A$. It is
convenient to set $x_{\emptyset ,i}=x_i$, $i=1,\dots ,n$.

Any ordering $i_1, i_2,\dots , i_n$ of indices $1,2,\dots , n$
defines a factorization of $P(t)$. Namely, let $A_k=\{i_1,
i_2,\dots , i_k\}$, $k=1,\dots , n-1$ and $A_0=\emptyset $. The
following result was obtained in \cite {GR2}.

\proclaim {Theorem 1.1.1}
$$P(t)=(t-x_{A_{n-1}, i_n})(t-x_{A_{n-2}, i_{n-1}})\dots (t-x_{A_0, i_1}).$$
\endproclaim

The elements $x_{A, i}$ (called {\it pseudo-roots} of $P(t)$ in
Section 1.2) are defined by pairs $A\subseteq \{1,2,\dots , n\}$
and $i\in \{1,2,\dots , n\}\setminus A$. They satisfy the
relations

$$x_{A\cup \{i\}, j}+x_{A, i}=x_{A\cup \{j\}, i}+x_{A, j}\ \ \ (1.1a)$$
$$x_{A\cup \{i\}, j}x_{A, i}=x_{A\cup \{j\}, i}x_{A, j}\ \ \ \ \ (1.1b)$$

When $A=\emptyset $, the relations (1.1) give relations (0.3).

To understand these relations better, we consider the Hasse graph
of the Boolean lattice of subsets of $\{1,2,\dots , n\}$. Thus, we
consider the directed graph whose vertices are subsets of
$\{1,2,\dots , n\}$ and an edge goes from vertex $B$ to vertex
$A$ if and only if $B=A\cup \{i\}$ where $i\in \{1,2,\dots ,
n\}\setminus A$. Then the relations (1.1) can be described nicely
by a {\it diamond} in the graph with vertices $A\cup \{i,j\},
A\cup \{i\}, A\cup \{j\}, A$ (see Fig 1 for $A=\emptyset $,
$i=1$, $j=2$.)

The elements $x_{A\cup \{i\}, j}, x_{A\cup \{j\}, i}, x_{A, j},
x_{A,i}$ correspond to the edges from $A\cup \{i,j\}$ to $A\cup
\{i\}$, from $A\cup \{i, j\}$ to $A\cup \{j\}$, from $A\cup
\{j\}$ to $A$, and from $A\cup \{i\}$ to $A$.

\subhead 1.2. Pseudo-roots\endsubhead The elements $x_{A, i}$
introduced in the previous section are natural examples of {\it
pseudo-roots} of noncommutative polynomials.

Let $R$ be an associative ring with unit and
$P(t)=a_0t^n+a_1t^{n-1}+\dots +a_n$ be a polynomial in $R[t]$.
According to \cite {GGRSW} an element $x\in R$ is {\it a
pseudo-root} of $P(t)$ if there exist polynomials $Q_1(t),
Q_2(t)\in R[t]$ such that
$$P(t)=Q_1(t)(t-x)Q_2(t).$$

The element $x$ is {\it a right root} of $P(t)$ if $Q_2(t)=1$ and
is {\it a left root} of $P(t)$ if $Q_1(t)=1$. As we mentioned in
the introduction, it is easy to check that $x$ is a right root of
$P(t)$ if and only if $a_0x^n+a_1x^{n-1}+\dots +a_n=0$.
Similarly, $x$ is a left root if and only if
$x^na_0+x^{n-1}a_1+\dots +a_n=0$.

For a noncommutative ring $R$, a theory of polynomials over $R$ should
be based not only on properties of right (left) roots but on
pseudo-roots as well (see Section 1.1.) In particular, it is
natural to study subrings  $R_P$ of $R$ generated by pseudo-roots
of a given polynomial $P(t)\in R[t]$. We will construct now a {\it
universal algebra of pseudo-roots} $Q_n$ introduced in \cite
{GRW}. In certain situations there exists a canonical
homomorphism of $Q_n$ into $R_P$.

\subhead 1.3. The algebra $Q_n$\endsubhead Let $F$ be a field. The algebra
$Q_n$ over $F$ is generated by elements $x_{A,i}$ where
$A\subseteq \{1,2,\dots , n\}$ and $i\in \{1,2,\dots ,
n\}\setminus A$ satisfying relations (1.1).
%$$z_{A\cup \{i\}, j}+z_{A, i}=z_{A\cup \{j\}, i}+z_{A, j},\ \ \ (2.2a)$$
%$$z_{A\cup \{i\}, j}z_{A, i}=z_{A\cup \{j\}, i}z_{A, j}\ \ \ (2.2b)$$

Let $R$ be an algebra over a field $F$ and let
$P(t)=t^n+a_1t^{n-1}+\dots +a_n$ be a polynomial over $R$ such
that $P(t)$ has $n$ right roots $x_1,\dots , x_n$ in a generic
position. The subalgebra of of $R$ generated by all pseudo-roots of $P(t)$ is
denoted by $R_P$. The following theorem was proved in \cite {GRW}.

\proclaim {Theorem 1.3.1} There exists a homomorphism
$$\alpha :Q_n\rightarrow R_P$$ such that $\alpha (z_{A,i})=x_{A,i}$.
\endproclaim

It turns out that the algebra $Q_n$ has several nice properties
studied in \cite {GRW, GGR, GGRSW, SW, Pi}. It is easy to see
that this algebra has a set of generators satisfying quadratic
relations because linear and quadratic relations in (1.1) are
separated, i.e. $Q_n$ is a quadratic algebra. In \cite {GRW} we
constructed a linear basis in $Q_n$, in \cite {GGRSW} we computed
the Hilbert series for $Q_n$ and its quadratic dual, and in \cite
{SW} we showed that $Q_n$ is a Koszul algebra.

Other properties of $Q_n$, including its relations to the algebra
of noncommutative symmetric functions introduced in \cite {GKLLRT},
were discussed in \cite {GRW}. Overall, one can say that $Q_n$ is
a rather ``tame" algebra (due to the properties listed above)
despite its exponential growth.

%Many properties of $Q_n$ can be explained by its connection with
%the Hasse graph of the Boolean lattice of subsets of the set
%$\{1,2,\dots , n\}$. We will show below that $Q_n$ is, in fact, a
%special case of a large class of algebras associated with
%directed graphs and that all mentioned properties of $Q_n$ can be
%generalized to this class of algebras. A ``geometric nature" of
%these algebras helps us to understand their structure and to
%simplify the proofs of main results compare to the proofs given
%in \cite {GRW, GGRSW, SW, Pi}.

\head 2. The algebras ${\bold A(i,\Gamma)}  $ and ${\bold
A(\Gamma)}$ \endhead

\subhead 2.1. Directed graphs\endsubhead Let  $\Gamma = (V, E)$ be
a {\bf directed graph}. That is, $V$ is a set (of vertices), $E$
is a set (of edges), and ${\bold t}: E \rightarrow V$ and ${\bold
h}: E \rightarrow V$ are functions. (${\bold t}(e)$ is the {\it
tail} of $e$ and ${\bold h}(e)$ is the {\it head} of $e$.)

A vertex $u\in V$ is called {\it maximal } if there is no $e\in E$
such that ${\bold h}(e)=u$. A vertex $v\in V$ is called {\it
minimal } if there is no $e\in E$ such that ${\bold t}(e)=v$.

We say that $\Gamma$ is {\bf layered} if $V$ is the disjoint union
of $V_i$, $0\leq i\leq n$, $E$ is the disjoint union of $E_i$,
$0\leq i\leq n$,
 ${\bold t}: E_i \rightarrow V_i$, \ ${\bold h}: E_i \rightarrow V_{i-1}$.
We will write $|v|=i$ if  $v\in V_i$. In this case the number $i$
is called the {\it level} of $v$. Note that a layered graph has no loops.
\smallskip
We will assume throughout the remainder of the paper that $\Gamma
= (V, E)$ is a layered graph with $V = \cup_{i=0}^n V_i$, and
$V_0 = \{*\}$ where $*$ is the unique minimal vertex of $\Gamma$.

To any partially ordered set $I$ there corresponds the directed
graph $\Gamma _I$ called the {\it Hasse graph} of $I$. Its
vertices are elements $x\in I$ and its edges are pairs $e=(x,
y)\in I\times I$ such that $y$ is an immediate predecessor of $x$
(in other words, $y<x$ and there is no $z\in I$ such that
$y<z<x$.) The element $x$ is the tail of $e$ and the element $y$
is the head of $e$.

Recall that a partially ordered set $I$ with a function
$r:I\rightarrow \Bbb Z_+$ is called a {\it ranked poset} with the
{\it ranking function} if $r(x)>r(y)$ for any $x>y$. Then $r$
turns the corresponding Hasse graph $\Gamma _I$ into a layered
graph with $|x|=r(x)$.
To any maximal (minimal) element in $I$ there corresponds a
maximal (minimal) vertex in $\Gamma _I$.

Note that every layered graph with no multiple edges arrives in 
in this manner. Although, the theory we develope below applies to
graphs with multiple edges (see Example 1 of Section 2.3), all
the examples we consider will be Hasse graphs.
\smallskip
Here are some examples of ranked partially ordered sets. 
\bigskip
\noindent {\bf Examples}.
\smallskip
 1. To any set $S$ corresponds the partially ordered set
$\Cal P(S)$ of all subsets of $S$. The order relation is given by
inclusion and $r(A)$ equals the cardinality of $A$. The minimal
element in $\Cal P(S)$ is the empty set and the maximal element is
$S$. We will always assume that $S$ is a finite set.

2. To any finite-dimensional vector space $E$ over a field $F$
corresponds the partially ordered set $\Cal W(E)$ of all vector
subspaces of $E$. The order relation is given by inclusion and
$r(W)$ is the dimension of $W$ for a subspace $W\subseteq E$.

The minimal element in $\Cal W(E)$ is the zero subspace $(0)$ and the
maximal element is $E$.

3. We say that a layered graph $\Gamma= (V,E)$ with $V = \cup
_{i=0}^nV_i$ is {\bf complete} if for every $i, 1 \le i \le n,$
and every $v \in V_i, w \in V_{i-1}$ , there is a unique edge $e$
with ${\bold t}(e) = v, {\bold h}(e) = w.$ A complete layered
graph is determined (up to isomorphism) by the cardinalities of
the $V_i$. We denote the complete layered graph with $V =
\cup_{i=0}^n V_i, |V_i| = m_i$ for $0 \le i \le n$, by $\bold
C[m_n,m_{n-1},\dots ,m_1,m_0]$.  Note that the graph $\bold
C[m_n,m_{n-1},\dots ,m_1,1]$ has a unique minimal vertex of level
$0$ and the graph $\bold C[1,m_{n-1},\dots ,m_1,m_0]$ has a unique
maximal vertex of level $n$.

4. Recall that the partially ordered set $\Cal Y$ of Young
diagrams can be identified with the set of weakly decreasing
sequences $\lambda =(\lambda _k)_{k\geq 1}$ such that $\lambda
_k=0$ for $k>>0$. By definition, $(\lambda _k)\geq (\mu _k)$ if
and only if $\lambda _k\geq \mu _k$ for all $k$. In fact, $\Cal
Y$ is a ranked partially ordered set. The rank is defined as
$r((\lambda _k))=\sum _{k\geq 1} \lambda _k$.

5. Abstract regular polytopes (see, for example, \cite {MS}) also
are natural examples of ranked partially ordered sets.

6. A family $\Cal F\subseteq \Cal P(S)$ is called {\it a complex}
if $B\in \Cal F$ and $A\subseteq B$ imply $A\in \Cal F$.  The
order and the ranking function on $\Cal P(S)$ induce an order and
a ranking function on $\Cal F$.
\smallskip
Another important example of a layered graph is the {\it graph of
right divisors of a monic polynomial} described in the next
subsection.

\subhead 2.2. The graph of right divisors \endsubhead Let $P(t)$
be a monic polynomial over an associative algebra $R$ and $S$ be a
set of pseudo-roots of $P(t)$. Denote by $R_S$ the subalgebra in
$R$ generated by pseudo-roots $x\in S$.

Construct a layered graph $\Gamma (P, S) =(V, E)$,
$$
V=V_n\cup V_{n-1}\cup \dots V_1\cup V_0
$$
as follows. The vertices of $V_k=\{v\in V: r(v)=k\}$ are monic
polynomials $B(t)\in R[t]$ such that $\deg B(t)=k$ and
$$
P(t)=Q(t)B(t)
$$
in $R[t]$.

We say that there is an edge from vertex $B_1(t)$ to $B_2(t)$ in
$\Gamma $ if
$$
B_1(t)=(t-x)B_2(t)
$$
for some $x\in S$.

Note that graph $\Gamma (P,S)$ has only one maximal vertex
$v=P(t)$ and only one minimal vertex $w=1$.

\subhead 2.3. Algebras ${\bold A(i,\Gamma)}  $ and ${\bold
A(\Gamma)}$ \endsubhead
\bigskip
In this section we discuss a class of algebras introduced in \cite {GRSW}
and \cite {RSW}.
Let  $\Gamma = (V, E)$ be a layered directed graph with a finite
number of layers. That is, $V = \cup _{i=0}^n V_i$, $E =
\cup_{i=1}^n E_i$,
 ${\bold t}: E_i \rightarrow V_i$, \ ${\bold h}: E_i \rightarrow V_{i-1}$.
\smallskip
Recall that  we are assuming throughout the remainder of the paper
that  $V_0 = \{*\}$ where $*$ is the unique minimal vertex of
$\Gamma$.

 For each $v \in \cup_{i=1}^n V_i$ we will fix, arbitrarily,
some $e_v \in E$, with ${\bold t}(e) = v$. Recall, that if $v \in
V_i$ we write $|v| = i$ and say that $v$ has {\it level} $i$.
Similarly, if $e \in E_i$ we write $|e| = i$ and say that $e$ has
{\it level} $i$.

If $v, \ w \in V$, a {\bf path} from $v$ to $w$ is a sequence of
edges $\pi = \{ e_1, e_2, ...,e_k \}$ with ${\bold t}(e_{1}) =
v$, ${\bold h}(e_k) = w$ and ${\bold t}(e_{i+1}) = {\bold
h}(e_i)$ for $1 \leq i < k$.  We write $v = {\bold t}(\pi)$, $w =
{\bold h}(\pi)$ and call $v$ the tail of the path and $w$ the
head of the path. We also write $v > w$ if there is a path from
$v$ to $w$.

With $\pi $ defined as above, let $l(\pi)=k$ be the {\it length}
of $\pi$ and let $|\pi|=|e_1| + ... + |e_k|$ be the {\it level} of
$\pi$.

If $\pi_1 = \{ e_1, ...,e_k \}$,$\pi_2 = \{ f_1, ...,f_l \}$ are paths with ${\bold h}(\pi_1) = {\bold t}(\pi_2)$
then $\{ e_1, ...,e_k, f_1, ..., f_l \}$ is a path; we denote it by $\pi_1\pi_2$.

For $v \in V$, write $v^{(0)} = v$ and define $v^{(i+1)} = {\bold h}(e_{v^{(i)}})$ for $0 \leq i <  |v|$.
Then $v^{(|v|)} = *$ and $\pi_v = \{ e_{v^{(0)}}, ...,e_{v^{(|v|-1)}} \}$  is a path from $v$ to $*$.

\smallskip Let $T(E)$ denote the free associative algebra  on $E$
over a field $F$. We are going to introduce a quotient algebra of
$T(E)$ modulo relations generalizing relations (0.2) and (0.3).
We will do this by equating coefficients of polynomials associated
with pairs of paths with the same origin and the same end.
%Define
%$$T(E)_i = span \{ e_1...e_r \ | \ r \geq 0, |e_1| + ... + |e_r| \leq i \}.$$

%If $a \in T(E)_i$, $a \not\in T(E)_{i-1}$, write $|a| = i$.
For a path $\pi = \{ e_1, e_2, ...,e_m \}$ define
$$P_{\pi}(\tau ) = (1-\tau e_1)...(1-\tau e_m) \in T(E)[\tau ]/(\tau ^{n+1}).$$

Note that $P_{\pi_1\pi_2}(\tau ) = P_{\pi_1}(\tau )P_{
\pi_2}(\tau )$ if ${\bold h}(\pi_1) = {\bold
t}(\pi_2)$. Write
$$P_{\pi}(\tau ) = \sum_{j=0}^{n+1} (-1)^j e(\pi,j)\tau ^j .$$
\smallskip \noindent
{\bf Definition 2.3.1.}  Let
$R$ be the ideal in $T(E)$ generated by
$$\{ e(\pi_1,k) - e(\pi_2,k) \ | \ k\geq 1, 
{\bold t}(\pi_1)={\bold t}(\pi_2),
{\bold h}(\pi_1)={\bold h}(\pi_2)\}$$

Set $$A(\Gamma) = T(E)/R.$$

In other words, $A(\Gamma )$ is defined by generators $e\in E$
and relations
$$\sum _{j=1}^ke_j=\sum _{j=1}^kf_j,$$
$$\sum _{i<j}e_ie_j=\sum _{i<j}f_if_j,$$
$$\dots $$
$$e_1e_2\dots e_k=f_1f_2\dots f_k.$$

\smallskip
\noindent {\bf Examples}. 1. Suppose $\Gamma = (V,E)$, $e,f\in E$
and $\bold t(e)=\bold t(f)$, $\bold h(e)=\bold h(f)$. Then the
images of $e$ and $f$ in $A(\Gamma )$ are equal. Thus 
$A(\Gamma )$ is isomorphic to $A(\Gamma ')$ where
$\Gamma '$ is the graph without multiple edges obtained from $\Gamma $
by identifying all edges with the same tail and head.
 
2. If $\Gamma $ is a tree-like graph,
i.e. there are no distinct paths $\pi _1, \pi _2$ such that
${\bold t}(\pi _1)={\bold t}(\pi _2)$,
${\bold h}(\pi _1)={\bold h}(\pi _2)$, then
$A(\Gamma )=T(E)$
is the free associative algebra generated by $e\in E$.

3. Let $\Gamma $ be a diamond graph from Fig 1. Then $A(\Gamma )$
is defined by the relations
$$e_1+e_2=f_1+f_2,$$
$$e_1e_2=f_1f_2$$
and $A(\Gamma )$ is isomorphic to the algebra $Q_2$ defined by the
relations (1.1).
\smallskip

It also useful to consider larger algebras $A(i,\Gamma )$ associated
to directed graphs. To introduce these classes of algebras, define 
$P_{i,\pi}(\tau )$ to be the image of $P_{\pi}(\tau )$ in the quotient
of $T(E)[\tau ]/(\tau ^{i+1})$ and write
$$P_{i,\pi}(\tau ) = \sum_{k=0}^{min(l(\pi),i)} (-1)^k e(i,\pi,k)\tau ^k $$
for $i\geq 1$.

Note that $P_{i, \pi_1\pi_2}(\tau ) = P_{i, \pi_1}(\tau )P_{i,
\pi_2}(\tau )$ if ${\bold h}(\pi_1) = {\bold t}(\pi_2)$. 

Set $e(i,\pi,k) = 0$ if $k > min(l(\pi),i)$. For $v \in
\cup_{l=1}^n V_l$, set $P_{i,v}(t) = P_{i,\pi_v}(t)$ and $e(i,v,k)
= e(i,\pi_v,k).$ Also, set $P_{i,*}(t) = 1$ and $e(i,*,k)= 0$ if
$k > 0.$
\smallskip \noindent
{\bf Definition 2.3.2.}  Let
$R(i)$ be the ideal in $T(E)$ generated by
$$\{ e(i,\pi_1,k) - e(i,\pi_2,k) \ | \ {\bold t}(\pi_1)={\bold t}(\pi_2), {\bold h}(\pi_1)=
{\bold h}(\pi_2), \ 1 \leq k \leq min(l(\pi_1),i) \}.$$
Note that
$$R(1) \subseteq R(2) \subseteq ... \subseteq  R(n) = R(n+1) = ...$$

Let $$A(i,\Gamma) = T(E)/R(i).$$
Therefore we have
$$A(1,\Gamma) \rightarrow A(2,\Gamma) \rightarrow ... \rightarrow A(n-1,\Gamma) \rightarrow 
A(n,\Gamma) .$$

Note that $A(n,\Gamma) = A(\Gamma)$ (as in Definition 2.3.1). 

In other words, $A(i,\Gamma )$ is defined by generators $e\in E$
and relations
$$\sum \Sb j_1<j_2<\dots <j_k\\ 1\leq k\leq i\endSb
(e_{j_1}e_{j_2}\dots e_{j_k}- f_{j_1}f_{j_2}\dots f_{j_k})$$ for
all pairs of paths $\pi _1=(e_1, e_2,\dots e_s)$, $\pi _2=(f_1,
f_2,\dots f_s)$ such that ${\bold t}(\pi _1)={\bold t}(\pi _2)$,
${\bold h}(\pi _1)={\bold h}(\pi _2)$.
\smallskip
 In Section 3.3 we will show that for most of the graphs
described in Section 2.1, algebras $A(2, \Gamma )$ and $A(\Gamma
)$ coincide. In this case $A(\Gamma )$ is described by linear and
quadratic relations. In particular, for the Hasse graph $\Gamma
_n$ of subsets of $\{1,2,\dots , n\}$ the algebra $A(\Gamma _n)$ is
isomorphic to the algebra $Q_n$ described in Section 1.3.

\subhead 2.4. Universality of $A(\Gamma )$\endsubhead In
notations of Section 2.3 assume that graph $\Gamma $ has a unique
maximal vertex $M$ and a unique minimal vertex $*$. Define a
polynomial $\Cal P_{\Gamma }(t)$ over algebra $A(\Gamma )$
corresponding to any path $\pi _0=(e_1,e_2,\dots, e_n)$ from $M$
to $*$:
$$\Cal P_{\Gamma }(t)=(t-e_1)(t-e_2)\dots (t-e_n).$$

It follows from the the definition of $A(\Gamma )$ that $\Cal
P_{\Gamma }(t)$ does not depend of a choice of path $\pi _0 $.

The polynomial $\Cal P_{\Gamma }(t)$ is a monic polynomial. We call it
the {\it universal polynomial} over $A(\Gamma )$.

 Let $R$ be an algebra, $P(t)$ a monic polynomial of
degree $n$ over $R$, and $S$ a set of pseudo-roots of $P(t)$. Let
$\Gamma (P,S)$ be the layered graph constructed in Section 2.2.

Assume that the set $S$ contains pseudo-roots $a_1, a_2,\dots ,
a_n$ such that $$P(t)=(t-a_1)(t-a_2)\dots (t-a_n).$$ Then  graph
$\Gamma (P,S)$ contains a directed path from maximal vertex
$M=P(t)$ to minimal vertex $*=1$.

Following Section 2.3 construct the algebra $A(\Gamma (P,S))$. Let
$\Cal P(t)=P_{\Gamma }(t)$ be the universal polynomial over this algebra.

\proclaim {Theorem 2.4.1} There is a canonical homomorphism
$$
\kappa : A(\Gamma (P, S))\rightarrow R
$$
such that the induced homomorphism of polynomial algebras
$$
\hat \kappa : A(\Gamma (P,S))[t]\rightarrow R[t]
$$
maps $\Cal P(t)$ to $P(t)$.
\endproclaim

Theorem 2.4.1 was proved in \cite {GGRW1}. To construct the
homomorphism $\kappa $, note that to any edge $e\in \Gamma (P,S)$
corresponds a pair of polynomials $B_1(t), B_2(t)$ in $R[t]$ such
that $B_1(t), B_2(t)$ divide $P(t)$ from the right and
$B_1(t)=(t-a)b_2(t)$. Set $\kappa (e)=a$. One can see that
 $\kappa $ can be uniquely extended to the homomorphism
$ A(\Gamma (P, S))\rightarrow R$ and that $\hat \kappa (\Cal
P)(t)=P(t)$.

\head 3. Properties of the algebras $A(i, \Gamma)$ and $A(\Gamma
)$\endhead

Throughout this section that we continue to assume $\Gamma
= (V, E)$ is a layered graph with $V = \cup_{i=0}^n V_i$, that
$V_0 = \{*\}$, and that, for every $v \in  V_+ = \cup_{i=1}^n
V_i$, $\{ e \in E \ | \ {\bold t}(e) = v \} \neq \emptyset$. For
each $v \in V_+$ fix, arbitrarily, some $e_v \in E$ with ${\bold
t}(e_v) = v$.
This defines a distinguished path, denoted by $\pi _v$, from $v$ to $*$.
Namely, for $v \in V^+$ we define $v^{(0)} =v$ and 
$v^{(i+1)} = {\bold h}(e_{v^{(i)}})$ for $0 \le i <|v|$
and we set $\pi _v=\{e_{v^{(0)}}, e_{v^{(1)}},\dots ,e_{v^{(|v|-1)}}\}$.  

\subhead 3.1. Linear basis in $A(\Gamma )$\endsubhead
 For $v \in V_+$ and $1 \le k \le |v|$ we define
$\hat{e}(v,k)$ to be the image in $A(\Gamma)$ of the product
$e_1\dots e_k$ in $T(E)$ where $\pi_v = \{e_1,\dots ,e_{|v|}\}$.

If $(v,k), (u,\ell) \in V \times {\Bbb N}$ we say $(v,k)$  {\bf
covers}  $(u,\ell)$ if $v > u$ and $k=|v| - |u|$. In this
case we write $(v,k) \gtrdot (u,\ell)$.  (In \cite{GRSW} we used
different terminology and notation: if $(v,k) \gtrdot (u,\ell)$ we
said $(v,k) $ can be  composed with $(u,\ell)$ and wrote $(v,k)
\models (u,\ell)$.)

\smallskip \noindent
{\bf Example}. In Fig 1  pair $(\{12\}, 1)$ covers the pairs
$(\{1\}, 1)$ and $(\{2\}, 1)$.
\smallskip
The following theorem is proved in \cite {GRSW} (see Corollary
4.5).

\proclaim{Theorem 3.1.1} Let $\Gamma = (V,E)$  be a layered graph,
$V = \cup_{i=0}^n V_i, $ and $V_0 = \{*\}$ where $*$ is the unique
minimal vertex of $\Gamma$. Then
$$\{\hat{e}(v_1,k_1)\dots \hat{e}(v_l,k_l)|l \ge 0, v_1,\dots ,v_l
\in V_+, 1 \le k_i \le |v_i|, (v_i,k_i) \not\gtrdot
(v_{i+1},k_{i+1})\}$$ is a basis for $A(\Gamma).$
\endproclaim

Theorem 3.1.1 implies the construction of a linear basis in
the algebra $Q_n$ obtained in \cite {RSW}, but the current description
is much nicer and the proof is much shorter.

 \subhead 3.2. New description of the algebras $A(i,\Gamma)$ and
$A(\Gamma )$\endsubhead In Section 2.3 we described algebras
$A(i,\Gamma )$ as algebras whose generators are edges of $\Gamma
$ subject to homogeneous relations of degree $1,2,\dots , i$. In
this Section, following \cite {RSW}, we describe those algebras as
algebras whose generators are vertices of $\Gamma $ subject to
homogeneous relations  of degree $2,\dots , i$. In fact, we will
represent each vertex $v$ in $\Gamma $ by a path from $v$ to $*$.
By Definition 2.3.1 the sum of edges in each of such paths has
the same image in $A(\Gamma )$.

Define $e(v,1) = e_{v^{(0)}} + e_{v^{(1)}} + ... +e_{v^{(|v|-1)}}.$ 
Thus
$$e_v = e(v,1) - e(v^{(1)},1) = e({\bold t}(e_v),1) - e({\bold h}(e_v),1).$$

Let $E' = \{e_v |v \in V^+\}.$ Let $F$ be our ground field. For
any set $X$ denote by $FX$ the vector space over $F$ with basis
elements $x\in X$. Define $\tau:FE \rightarrow FE'$ by
$$\tau(f) = e({\bold t}(f),1) - e({\bold h}(t),1).$$
Then $\tau$ is a projection of $FE$ onto $FE'$ with kernel
generated by $e(i,\pi _1, 1)-e(i,\pi _2, 1)$ where $\pi _1, \pi
_2$ are paths with the same tail and head.

Now define $\eta: FE' \rightarrow FV^+$ by
$$\eta: e_v \mapsto v.$$

Then $\eta$ is an isomorphism of vector spaces and $\eta\tau$
induces a surjective homomorphism of graded algebras
$$\theta: T(E) \rightarrow T(V^+).$$

\proclaim {Proposition 3.2.1} $\theta $ induces an
isomorphism $A(i, \Gamma) \cong T(V^+)/\theta(R(i)).$
\endproclaim

It is important to write generators for the ideal $\theta(R(i))$
explicitly. The ideal is generated by elements of the form
$\theta (e(i,\pi_1,k) - e(i,\pi_2,k))$ where $2\leq k\leq i$ (see
Section 2.3 for the definition of $e(i,\pi_1,k))$. Therefore it
will be sufficient to write a formula for $\theta(e(i,\pi,k))$.
Let $\pi = \{e_1,e_2,...,e_m\}$ be a path, let ${\bold{t}}(e_j) =
v_{j-1}$ for $1 \le j \le m$ and let ${\bold h}(e_m ) = v_m.$

\proclaim{Proposition 3.2.2}
$$\theta(e(i,\pi,k)) = (-1)^k\sum_{1 \le j_1 < ... < j_k \le s} (v_{j_1 - 1}- v_{j_1})...(v_{j_k - 1 }-v_{j_k})$$
where $s=\min (i,m)$.
\endproclaim

Proposition 3.2.1 immediately implies
 \proclaim{Corollary 3.2.3} $A(2,\Gamma ) $ is a quadratic
algebra.
\endproclaim

We will show below that for many interesting graphs algebra
$A(2,\Gamma )$ coincides with algebra $A(\Gamma )$.

\subhead 3.3. Algebras associated with uniform graphs
\endsubhead

Let $\Gamma$ be a layered graph. For $v \in V$ define $\Cal
S_-(v)$ to be the set of all vertices  $w \in V$ covered by $v$,
i.e. such that there exists an edge with the tail $v$ and the
head $w$. Similarly, define $\Cal S_+(v)$ to be the set of all
vertices covering $v$.

For $v \in V_j, j \ge 2$, let $\sim_v$ denote the equivalence
relation on ${\Cal S}_-(v)$ generated by $u \sim_v w $ if ${\Cal
S}_-(u) \cap {\Cal S}_-(w) \ne \emptyset.$

\definition{ Definition 3.3.1} The layered graph $V$ is said to be {\bf uniform} if, 
for every $v \in V_j, j \ge 2$,
all elements  of ${\Cal S}_-(v)$ are equivalent under $\sim_v$.
\enddefinition

\smallskip \noindent
{\bf Example}. The diamond graph in Fig 1 is a uniform graph.
\smallskip

 \proclaim{Proposition 3.3.2} The Hasse graphs of
partially ordered sets listed in Examples 1-5 in
Section 2.1 are uniform.
\endproclaim

Algebras associated to uniform graphs have especially nice
structure.

\proclaim{Proposition 3.3.3} Let $\Gamma$ be a uniform layered
graph.  Then $A(\Gamma) =A(2,\Gamma )\cong T(V^+)/R_V$ is a
quadratic algebra and $R_V$ is generated by
$$\{v(u-w) - u^2 + w^2 + (u-w)x| v \in \cup_{i=2}^n V_i, u,w \in
{\Cal S}_-(v), x \in {\Cal S}_-(u) \cap {\Cal S}_-(w)\}.$$
\endproclaim

\proclaim{Corollary 3.3.4} The algebra associated to the Hasse graph
of all subsets of $\{1,2,\dots, n\}$ is isomorphic to algebra
$Q_n$ described in Section 1.3
\endproclaim

\subhead 3.4. $A(\Gamma)$ is a Koszul algebra
\endsubhead
Koszul algebras constitute an important class of quadratic
algebras (see, for example,  \cite {PP, U} for the definitions
and properties of Koszul algebras). Most of the known examples are
algebras of a polynomial growth. A wide variety of Koszul algebras
of exponential growth is given by the following theorem proved in
\cite {RSW}.

\proclaim {Theorem 3.4.1} Let $\Gamma$ be a uniform layered graph
with a unique minimal element. Then $A(\Gamma)$ is a Koszul
algebra.
\endproclaim

\subhead 3.5. Algebras dual to $A(2,\Gamma )$\endsubhead

Let $A$ be a quadratic algebra over a field $F$. Thus
$A$ is isomorphic to a quotient algebra $T(W)/<L>$ where $W$ is a
vector space over $F$ and $L$ is a subspace in $W\otimes W$.
Assume that $W$ is finite dimensional.
Denote by $W^*$ the dual vector space and  set $L^{\perp }=\{f\in
W^*\otimes W^*| f|_L=0\}$. By definition, the algebra
$A^!=T(W^*)/<L^{\perp }>$ is the quadratic dual to algebra $A$.

Both algebras $A$ and $A^!$ are graded and their Hilbert series
$H(A,\tau )$ and $H(A^!, \tau )$ are well defined. A graded algebra $A$ is
Koszul if and only if $A^!$ is Koszul and in this case the
Hilbert series of $A$ and $A^!$ are related by
$$ H(A,\tau )H(A^!, -\tau )=1$$
(see, for example, \cite {PP, U}.)

According to Corollary 3.2.3 the algebra $A(2, \Gamma )$ is quadratic.
Denote its quadratic dual algebra by $B(\Gamma )$. If the graph
$\Gamma $ is uniform than $A(2,\Gamma)=A(\Gamma)$ and $B(\Gamma
)=A(\Gamma )^!$.

The following theorem gives a description of $B(\Gamma )$. When
$\Gamma $ is the Hasse graph of the partially ordered set of subsets
of $\{1,2,\dots , n\}$ such description was essentially given in
\cite {GGRSW}.

Let $\Gamma =(V, E)$ be a layered graph with a unique minimal
vertex. Assume that $V$ is a finite set.

\proclaim{Theorem 3.5.1} The algebra $B(\Gamma )$ is generated by
vertices $v\in V$ subject to the following relations:

i) $uv=0$ if $u\neq v$ and there is no edge $e\in E$ such ${\bold
t}(e)=u$,
 ${\bold h}(e)=v$;

ii) $v^2+v\sum  w=0$, where the sum is taken over all $w\in \Cal
S_-(v)$;

iii) $v^2+(\sum u)v=0$,  where the sum is taken over all $u\in \Cal
S_+(v)$.
\endproclaim

In Section 4 we will compute Hilbert series for some algebras
$B(\Gamma )$ and see that those Hilbert series are polynomial in
$\tau $. It follows that in this case, the  algebras $B(\Gamma )$ are
finite dimensional and the ideal (of codimension one) generated by $V$ is 
nilpotent. It is not surprising since Theorem 3.5.1 implies that 
$v^3=0$ for any $v\in V$.

\head 4. The Hilbert series of $A(\Gamma)$\endhead

In this section  we compute the Hilbert series of algebras
$A(\Gamma)$ introduced in Section 2.3 and specialize the result
for some examples of layered directed graphs. All results
formulated in this section were obtained in \cite {RSW1}.

\subhead 4.1. Main theorem\endsubhead

Let $h(\tau )$  denote the Hilbert series $H(A(\Gamma),\tau )$,
where $\Gamma =(V,E)$ is a layered  graph with unique minimal element
$*$ of level $0$. Arrange the elements of $V$ in nonincreasing
order and index the elements of vectors and matrices by this
ordered set.
Let ${\bold 1}$ denote the column vector all of whose entries are $1$, and
let $\zeta(\tau )$ denote the matrix with entries
$\zeta_{v,w}(\tau )$ for $ v,w \in V$ where $\zeta_{v,w}(\tau )
= \tau ^{|v|-|w|}$ if $v \ge w$ and $0$ otherwise. Then we have

\proclaim{Proposition 4.1.1}
$$h(\tau ) = \frac{1-\tau }{1 - \tau {\bold 1}^T\zeta(\tau )^{-1}{\bold 1}}.$$
\endproclaim

Note that $\zeta (\tau )=1+N$ where $N$ is strongly upper-triangular.
Consequently,
the $(v,w)$-entry of $\zeta (\tau )^{-1}$ can be written as
$$\sum _{v=v_1>\dots >v_{\ell}=w\geq *}(-1)^{\ell +1}\tau ^{|v|-|w|}$$
and we have the following result.

\proclaim{Theorem 4.1.2} Let $\Gamma$ be a layered graph with
unique minimal element $*$ of level $0$ and $h(\tau ) $ denote the
Hilbert series of $A(\Gamma)$.  Then

$$h(\tau ) = 
\frac{1-\tau }{1+\sum_{v_1 > v_2 \dots > v_{\ell}\geq *} (-1)^{\ell}
\tau ^{|v_1| - |v_{\ell}|+1}}.$$
\endproclaim

The proof of Proposition 4.1.1 and Theorem 4.1.2 is based on Theorem 3.1.1 
describing a linear basis in $A(\Gamma )$.

We remark that the matrices $\zeta(1)$ and $\zeta(1)^{-1}$ are
well-known as the zeta-matrix and the
 M\"obius-matrix of $V$ (cf. \cite R).

In the remaining part of this section we will use Theorem 4.1.2 
to compute the Hilbert series of the algebras
$A(\Gamma)$ associated with certain layered graphs.

\subhead 4.2. The Hilbert series of the algebra associated with
the Hasse graph of the lattice of subsets of
$\{1,\dots,n\}$\endsubhead

Let $\Gamma_n$ denote the Hasse graph of the lattice of all
subsets of $\{1,\dots,n\}$.  Thus the vertices of $\Gamma_n$ are
subsets of $\{1,\dots,n\}$, the order relation $>$ is set
inclusion $\supset$, the level $|v|$ of a set $v$ is its
cardinality, and the unique minimal vertex $*$ is the empty set
$\emptyset$. Then the algebra $A(\Gamma_n)$ is the algebra $Q_n$
defined in \cite{GRW}. Theorem 4.1.2 implies the following theorem
(from \cite{GGRSW}). The proof obtained in \cite {RSW} is much
shorter and more conceptual than that in \cite{GGRSW}.

\proclaim{Theorem 4.2.1} $$H(Q_n,\tau ) =  \frac{1-\tau }{1 - \tau
(2-\tau )^n}.$$
\endproclaim

\subhead 4.3. The Hilbert series of  algebras associated with the
Hasse graph of the lattice of subspaces of a finite-dimensional
vector space over a finite field\endsubhead

We will denote by $\bold L(n,q)$ the Hasse graph of the lattice
of subspaces of an $n$-dimensional space over the field $\bold
F_q$ of $q$ elements.  Thus the vertices of $\bold L(n,q)$ are
subspaces of $\bold F_q^n$, the order relation  $>$ is inclusion
of subspaces $\supset$, the level $|U|$ of a subspace $U$ is its
dimension, and the unique minimal vertex $*$ is the zero subspace
$(0)$. Recall that ${\binom nm}_q$ is a $q$-binomial coefficient.

\proclaim{Theorem 4.3.1}
$$
\frac{1 - \tau }{H(A(\bold L(n,q)),\tau )} = 1 - \tau \sum_{m = 0}^n
{\binom nm}_q (1 - \tau )(1 - \tau q) \dots (1 - \tau q^{n - m -
1}).$$
\endproclaim

Note that setting $q=1$ in the expression in Theorem 4.3.1 gives
$1-\tau (2-\tau)^n$. By Theorem 4.2.1, this is $(1-\tau )H(Q_n, \tau
)^{-1}$.

 Since by
Theorem 3.4.1 (see also \cite{RSW}) $A(\bold L(n,q))$ is a Koszul
algebra, we have the following corollary.

\proclaim{Corollary 4.3.2} $$H(A(\bold L(n,q))^!,\tau ) = 1 +
\sum_{m=0}^{n-1} {\binom nm_q }(1 + \tau q)\dots(1+\tau
q^{n-m-1}).$$
\endproclaim
\smallskip
\subhead 4.4. The Hilbert series of algebras associated with
complete layered graphs\endsubhead

Recall (see Example 3, Section 2.1) that a layered graph $\Gamma=
(V,E)$ with $V = \cup _{i=0}^nV_i$ is {\bf complete} if for every
$i, 1 \le i \le n,$ and every $v \in V_i, w \in V_{i-1}$ , there
is a unique edge $e$ with ${\bold t}(e) = v, {\bold h}(e) = w.$ A
complete layered graph is determined (up to isomorphism) by the
cardinalities of the $V_i$. We denote the complete layered graph
with $V = \cup_{i=0}^n V_i, |V_i| = m_i$ for $0 \le i \le n$, by
$\bold C[m_n,m_{n-1},\dots ,m_1,m_0]$.  Note that the graph $\bold
C[m_n,m_{n-1},\dots ,m_1,1]$ has a unique minimal vertex of level
$0$ and so Theorem 4.1.2 applies to $A(\bold
C[m_n,m_{n-1},\dots ,m_1,1])$. This leads to the following
theorem.

\proclaim{Theorem 4.4.1}
$$\frac{1-\tau }{H(A(\bold C[m_n,m_{n-1},\dots ,m_1,1],\tau )} = $$
$$ 1 - \sum_{k = 0}^n  \sum_{a = k}^n  \; (-1)^k
m_a (m_{a-1} - 1)(m_{a-2} - 1) \dots (m_{a-k+1} - 1)m_{a-k} \tau
^{k+1}.$$
\endproclaim

When $k=1$ and $k=2$, the product 
$(m_{a-1} - 1)(m_{a-2} - 1) \dots (m_{a-k+1} - 1)$
in this expression (and also in the expression in Corollary 4.4.2) 
represents the empty product and so has value $1$.

This result applies, in particular, to the case $m_0 = m_1 = ...
= m_n = 1.$ The resulting algebra $A(\bold C[1,\dots,1])$ has $n$
generators and no relations. Theorem 4.4.1 shows that
$$\frac{1-\tau }{H(A(\bold C[1,\dots,1]),\tau)}
= 1 - \sum_{a=0}^n \tau + \sum_{a=1}^n \tau ^2 = (1-\tau
)(1-n\tau).$$

Thus $$H(A(\bold C[1,\dots,1]),\tau ) = \frac1{1-n\tau }$$ and we
have recovered the well-known expression for the Hilbert series
of the free associative algebra on $n$ generators.

Since by Theorem 3.4.1 (see also \cite{RSW})  the algebras
associated to complete directed graphs are Koszul algebras,  we
have the following corollary.

\proclaim{Corollary 4.4.2}
$$H(A(\bold C[m_n,m_{n-1},\dots ,m_1,1])^!,\tau )=
$$ $$1 +     \sum_{k=1}^n\sum_{a =
k}^n  \; m_a (m_{a-1} - 1)(m_{a-2} - 1) \dots (m_{a-k+1} - 1)\;
\tau ^k.$$
\endproclaim

\head 5. Sufficient sets of pseudo-roots and directed
graphs\endhead

We return to questions 1) and 2) from the introduction: Given a
polynomial $P(t)$ over a noncommutative algebra, how to obtain its
factorizations and how to relate two different factorizations? To
answer these questions we will work in general context of
algebras $A(\Gamma )$.

\subhead 5.1. Defining sets of pseudo-roots\endsubhead In this
section we briefly describe some results from \cite {GGRW1}. Let
$R$ be an associative ring with unit, $P(t)=t^n+a_1t^{n-1}+\dots
+a_n$ be a polynomial over $R$, and $t$ be a central variable.

Recall (see Theorem 1.1.1) that if $P(t)$ has right roots $x_1,
x_2,\dots , x_n$ in a generic position then $P(t)$ admits
factorizations
$$P(t)=(t-x_{A_{n-1}, i_n})(t-x_{A_{n-2}, i_{n-1}})
\dots (t-x_{\emptyset ,i_1})$$ indexed by orderings of
$\{1,2,\dots , n\}$. Also recall that pseudo-roots $x_{A,i}$'s
are rational expressions in $x_1, x_2,\dots , x_n$.

According to \cite {GRW} the pseudo-roots $x_{A,i}$'s can be obtained
from a generic set of $n$ right roots by a sequence of operations
$$a,b \mapsto (a-b)a(a-b)^{-1}, \ \ (b-a)b(b-a)^{-1} \tag 5.1a$$ 
and from a generic set of $n$ left
roots by a sequence of operations
$$a,b \mapsto (a-b)^{-1}a(a-b), \ \
(b-a)^{-1}b(b-a). \tag 5.1b$$

This leads us to the following natural question. We call a set of
pseudo-roots $Y=\{y_1, y_2,\dots , y_n\}$ of $P(t)$ is a {\it
defining set} if $$P(t)=(t-y_n)(t-y_{n-1})\dots (t-y_1).$$

Our question is then: given a set of pseudo-roots $Z$, when it is
possible to construct a
defining set of pseudo-roots from elements of
$Z$ by a successive application of operations of (5.1)-type?

Our answer to this question is based on a geometrical ``diamond"
interpretation of operations of (5.1)-type.

 \subhead 5.2. Sufficient sets of pseudo-roots and
the algebra $Q_n$\endsubhead

To avoid taking inverses, we will slightly change the definition
of operations of (5.1)-type. Recall that algebra $Q_n$ corresponds
to the Hasse graph $\Gamma _n$ of the Boolean lattice of
$\{1,2,\dots ,\}$ and $Q_n[t]$  contains a  unique universal polynomial,
denoted $\Cal P(t)$ (as defined in Section 2.4.)

\definition{Definition 5.2.1} We say that a pseudo-root
$\xi \in Q_n$ is obtained from an ordered pair of pseudo-roots
$x_{A,i}, x_{B,j}$ by the $u$-{\it operation} if $A\cup
\{i\}=B\cup \{j\}$ and $(x_{A,i}- x_{B,j})x_{A,i}= \xi
(x_{A,i}-x_{B,j})$.

A pseudo-root $\eta \in Q_n$ is obtained from an ordered pair of
pseudo-roots $x_{A,i}, x_{B,j}$ by the $d$-{\it operation} if
$A\cup \{i\}=B\cup \{j\}$ have and 
$(x_{A,i}- x_{B,j})\eta=x_{A,i}(x_{A,i}-x_{B,j})$.
\enddefinition

\proclaim {Proposition 5.2.2} The element $x_{A\cup \{i\}, j}$ is
obtained by the $u$-operation from the pair $x_{A, j}, x_{A,i}$.

The element $x_{A, i}$ is obtained by the $d$-operation from the
pair $x_{A\cup \{j\}, i}, x_{A\cup \{i\}, j}$.
\endproclaim

\definition
{Definition 5.2.3} The set of elements in $Q_n$ that can be
obtained from elements of $Z$ by a successive applications of
$d$- and $u$-operations is called the $du$-{\it envelope} of $Z$.
\enddefinition

\definition {Definition 5.2.4} A set $Z\subseteq Q_n$ is called sufficient
if the $du$-envelope of $Z$ contains a defining set of
pseudo-roots of $\Cal P(t)$.
\enddefinition

Any defining set of elements is a sufficient set. Other examples
of sufficient sets in $Q_n$ are given by the following statement.

\proclaim {Proposition 5.2.5} The sets $\{x_{\emptyset, k}\ |
1\leq k\leq n\}$ and
$\{x_{\{12\dots \hat k \dots n\}, k}\ |1\leq k\leq n\}$ are
sufficient in $Q_n$ for $\Cal P (t)$.
\endproclaim

A necessary condition for a subset in $Q_n$ to be sufficient for
$\Cal P(t)$ is given by the following theorem.

\proclaim {Theorem 5.2.6} If $Z=\{x_{A_1, i_1}, x_{A_2, i_2},
\dots , x_{A_n, i_n}\}$ is a sufficient subset of $Q_n$ then 
$i_1, i_2,\dots , i_n$ are distinct.
\endproclaim

A set of edges in a directed graph is connected if it is connected
in the associated non-directed graph (see Section 5.3 below for
details).

\proclaim {Theorem 5.2.7} Let $Z=\{x_{A_1, i_1},
x_{A_2, i_2}, \dots , x_{A_n, i_n}\}$ be a subset of $Q_n$ such
that $i_1$, $i_2$, $\dots $, $i_n$ are distinct. If the set of edges
$\{(A_1, i_1), (A_2, i_2), \dots , (A_n, i_n)\}$ in $\Gamma _n$
is connected then the set $Z$ is sufficient.
\endproclaim

 Let $f:Q_n \rightarrow D$ be a
homomorphism of $Q_n$ into a division ring $D$, $\hat f:Q_n[t]
\rightarrow D[t]$ be the induced homomorphism of the polynomial
rings and $P(t)=\hat f(\Cal P(t))$.

\proclaim {Corollary 5.2.8} Let $Z=\{x_{A_1, i_1}, x_{A_2, i_2},
\dots , x_{A_n, i_n}\}$ be a subset of $Q_n$ such that $i_1$,
$i_2$, $\dots $, $i_n$ are distinct and the set of edges $\{(A_1, i_1),
(A_2, i_2), \dots , (A_n, i_n)\}$ in $\Gamma _n$ is connected.
Then all coefficients of $P(t)\in D[t]$ can be obtained from
elements $f(z)$, $z\in Z$, by operations of addition, subtractions,
multiplication, and left and right conjugation.
\endproclaim

\example {Examples} 1. The set $X=\{x_{\emptyset , 1},
x_{\emptyset , 2}, \dots , x_{\emptyset , n}\}$ is connected (the
corresponding edges have a common head $\emptyset $). Therefore,
$X$ is a sufficient set.

2. For $n=2$ the sufficient sets are $\{x_{\{i\},j}, x_{\emptyset
,i}\}$, $\{x_{\emptyset ,j}, x_{\emptyset ,i}\}$, $\{x_{\{i\},j},
x_{\{j\},i}\}$. The sets $\{x_{\{i\},j}, x_{\emptyset ,j}\}$ are
not sufficient. Here $i,j=1,2$, $i\ne j$.

3. Let $n=3$. The set $\{x_{\{1\},2}, x_{\{2\},1}, x_{\{1\},3}\}$
is sufficient because
$$
(x_{\{2\},1}-x_{\{1\},2})x_{\emptyset,1}
=x_{\{1\},2}(x_{\{2\},1}-x_{\{1\},2}),
$$
$$
x_{\{12\},3}(x_{\{1\},3}-x_{\{1\},2})
=(x_{\{1\},3}-x_{\{1\},2})x_{\{1\},3}
$$
and $\{x_{\{12\},3}, x_{\{1\},2}, x_{\emptyset , 1}\}$ is a
defining set of pseudo-roots.

The sets $\{x_{\{1\},2}, x_{\{2\},1}, x_{\emptyset, 3}\}$ and
$\{x_{\{1\},3}, x_{\{1\},2}, x_{\emptyset , 1}\}$ also are
sufficient but not defining sets in $Q_3$.

4. The set $W=\{x_{\{12\},3}, x_{\{3\},2}, x_{\emptyset ,1}\}$ is
not sufficient because $d$- and $u$-operations are not defined on
elements of $W$.
\endexample

Theorem 5.2.7 follows from a more general theorem for algebras
associated with directed graphs.

\subhead 5.3. Sufficient sets of edges for directed
graphs\endsubhead

In this section we will define and study sufficient sets of edges
in directed graphs $\Gamma =(V,E)$. These sets will provide us
with a construction of sufficient sets of pseudo-roots of
polynomials $\Cal P(t)$ over algebras $A(\Gamma )$. All graphs
considered in this sections are {\it simple} (i.e., if
$t(e)=t(f)$ and $h(e)=h(f)$ then $e=f$) and {\it acyclic} (i.e.,
there are no directed paths $P$ such that $t(P)=h(P)$).

Let $\Gamma =(V, E)$ be a directed graph.

\definition {Definition 5.3.1}
\roster
\item A pair of edges  $f_1, f_2$ with a common head
is obtained from the pair  $e_1, e_2$ with a common tail by
$D$-operation if $h(e_i)=t(f_i)$ for $i=1,2$;

\item A pair of edges  $e_1', e_2'$ with a common tail
is obtained from the pair $f_1', f_2'$ with a common head by
$U$-operation if $h(e_i')=t(f_i')$ for $i=1,2$.
\endroster
\enddefinition
\smallskip \noindent
{\bf Example}. In Fig 1 edges $f_1, f_2$ can be obtained from
edges $e_1, e_2$ by a $D$-operation. Conversely, edges $e_1, e_2$
can be obtained from edges $f_1, f_2$ by an $U$-operation.
\smallskip
\smallskip \noindent
{\bf Remark}. We do not require the uniqueness of $D$- and
$U$-operations.
\smallskip

\definition{Definition 5.3.2} A subset $E_0\subseteq E$ is called
$DU$-complete (or simply complete) if the results of any
$D$-operation or any $U$-operation applied to edges from $E_0$
belong to $E_0$.
\enddefinition

\proclaim {Proposition 5.3.3} For any subset $F\subseteq E$ there
exists a {\it minimal} $DU$-complete set $\hat F\subseteq E$
containing $F$.
\endproclaim

We call $\hat F$ the {\it completion} of $F$.

 Let $\Gamma =(V, E)$ be a directed graph.

\definition {Definition 5.3.4} A set of edges $G$ in $\Gamma $ is called
{\it sufficient\/} if its completion $\hat G$ contains a path
from a maximal vertex (source) to a minimal vertex (sink).
\enddefinition

\definition
{Definition 5.3.5} A set of vertices $W\subseteq V$ is called {\it
ample} if \roster
\item For any non-minimal vertex $v\in V$ there exists a vertex
$u\in W$ such that there is no directed path in $\Gamma $ from
$u$ to $v$;
\item
For any non-maximal vertex $v\in V$ there exists a vertex $w\in W$
such that there is no directed path in $\Gamma $ from $v$ to
$w$.
\endroster
\enddefinition

A set of edges is called ample if the set of its tails and heads
is ample.

\smallskip

As an example, consider the graph $\Gamma_n$ of all subsets of
$\{1,\dots,n\}$. It has one source $\{1,\dots,n\}$ and one sink
$\emptyset$.

\proclaim {Proposition 5.3.6}  A set of edges  $(A_1, i_1), (A_2,
i_2), \dots , (A_n, i_n)$ in $\Gamma_n $ is ample if $i_1$,
$i_2$, $\dots $, $i_n$ are distinct.
\endproclaim

\definition {Definition 5.3.7} A directed graph is called a modular graph if:

\roster
\item  For any two edges $e_1, e_2$ with a common tail there exist
edges $f_1, f_2$ with a common head such that $h(e_i)=t(f_i)$ for
$i=1,2$;
\item  For any two edges $h_1, h_2$ with a common head there exist
edges $g_1, g_2$ with a common tail such that $h(g_i)=t(h_i)$ for
$i=1,2$.
\endroster
\enddefinition
We do not require the uniqueness of $f_1, f_2$ and $h_1, h_2$.

\proclaim {Theorem 5.3.8} Any ample connected set of edges  of a
finite modular directed graph is a sufficient set.
\endproclaim

Now let $\Gamma =(V, E)$ be a directed graph such that \roster
\item  $\Gamma $ contains a unique source $M$ and a unique
sink $m$;
\item For each vertex $v\in V$ there exist a directed path from $M$ 
to $v$ and a directed path from $v$ to $M$.
\endroster
\smallskip
Recall that we associate to $\Gamma $ an algebra $A(\Gamma )$ and
the universal polynomial $\Cal P(t)\in A(\Gamma )[t]$. The polynomial $\Cal
P(t)$ is constructed using a path $e_1, e_2, \dots , e_n$ from
$M$ to $m$ in $\Gamma$, but it does not depend on the path. To
any edge $e\in E$ there corresponds to a pseudo-root $e\in
A(\Gamma )$ of $\Cal P(t)$, and to any path $(e_1, e_2,\dots ,
e_n)$ from $M$ to $m$ in $\Gamma$ there corresponds the
factorization
$$\Cal P(t)=(t-e_1)(t-e_2)\dots (t-e_n)$$
of $\Cal P(t)$ over $A(\Gamma )$.

Theorem 5.3.8 implies

\proclaim {Theorem 5.3.9} Let $S\subseteq E$ be an ample connected
set of edges in a modular directed graph $\Gamma =(V,E)$. Then
there exists a factorization (3.1) of $\Cal P(t)$ such that
the $DU$-completion of $S$ contains elements $e_1, e_2,\dots , e_n$
and, therefore, coefficients of $\Cal P(t)$.
\endproclaim

\enddocument

\bye